\documentclass[11pt]{article}
\usepackage{epsfig}
\newlength{\textbildsep}
\usepackage{amssymb}
\usepackage{amsmath}
\usepackage{algorithm}
 \usepackage{algpseudocode}

\usepackage{hyperref}

\newtheorem{theorem}{Theorem} 
\newtheorem{lemma}{Lemma}
\newtheorem{corollary}{Corollary} 
\newtheorem{remark}{Remark}
\newtheorem{proposition}{Proposition}

\newtheorem{definition}{Definition}
\newtheorem{fact}{Fact}
\newtheorem{conjecture}{Conjecture}

\input mssymb.tex

\begin{document}
\title{Tunnels in digital objects}
\author{Valentin E. Brimkov\thanks{Mathematics Department, SUNY Buffalo State University,
1300 Elmwood Ave., Buffalo, NY 14222.\\
Email: brimkove@buffalostate.edu.
}}


\title{On Guarding Line Arrangements
}






\date{}
\maketitle
\begin{abstract}
We study the problem of guarding a set of straight lines in the plane by selecting a minimum-cardinality set of points such that every line contains at least one selected point. 
This problem is known to be strongly NP-hard \cite{megiddo}.
We investigate approximation algorithms for the problem and analyze their performance guarantees 
under various structural assumptions on the underlying line arrangement. In particular, we present 
approximation schemes based on set-cover techniques, relaxation methods, and properties of maximal lines. 
We establish constant-factor bounds on the approximation ratio for several classes of arrangements and derive relationships between the guarding number and the incidence parameters of the arrangement. Connections with Dirac-type incidence results are also explored.
\smallskip


\end{abstract}

\section{Introduction}

Although straight lines are among the simplest geometric objects, several
challenging problems defined on sets of lines continue to attract attention
in discrete geometry (see, e.g., \cite{grunbaum1} and Chapter~7.2 of
\cite{preparata}). Some of these problems have resisted resolution for
decades; a notable example is the Dirac--Motzkin conjecture, recently
settled for large sets of lines by Green and Tao \cite{tao}. Other
questions remain open.

The present work concerns several problems related to line arrangements.
One of the central problems considered in this paper is the following:

\begin{description}
\item[Guarding a Set of Lines (GSL).]
\emph{Given a set of straight lines in the plane, not all of them parallel,
find a set of points of minimum cardinality such that each line contains
at least one of these points.}
\end{description}

The GSL problem can be interpreted as determining the minimum number of
locations for guards, cameras, or sensors required to monitor a network
of streets, corridors, tunnels, or pipes. Thus, from the standpoint of
applications, the problem is related to surveillance and security.

Let $g(M)$ denote the minimum number of points (or ``guards'') that cover
all lines in a set $M$. Megiddo and Tamir showed that GSL is strongly
NP-hard \cite{megiddo}; therefore, computing $g(M)$ is strongly NP-hard
as well. A natural objective is therefore to design approximation
algorithms for GSL, which is one of the goals of this paper.

Many geometric hitting set and set cover problems admit constant-factor approximation algorithms under additional structural assumptions. 
Typical approaches rely on properties such as bounded VC-dimension \cite{BronnimannGoodrich}, the existence of small $\varepsilon$-nets \cite{HausslerWelzl}, 
bounded union (or shallow-cell) complexity \cite{ClarksonVaradarajan}, or bounded overlap density. Although these approaches are based on different 
geometric and combinatorial characteristics, they all exploit restrictions on the complexity of the associated range space that are absent in general instances of Set Cover.

The approach developed in this paper is of a different nature.
Rather than relying on global complexity measures of the associated
range space, our analysis exploits incidence properties of line
arrangements. In particular, the approximation guarantees are derived
from bounds involving the degrees of intersection points and the
frequencies of lines, that is, the numbers of intersection points
lying on a line. Consequently, the approximation guarantees are
expressed directly in terms of incidence parameters of the
arrangement, most notably the degrees of intersection points and the
maximum line frequency $w(M)$.

Given a set of lines $M$, not all parallel, let $w(M)$ denote the maximum
number of intersection points lying on any single line of $M$. One of the
questions we investigate concerns the relationship between $g(M)$ and
$w(M)$. The quantity $w(M)$ is closely related to the weak dual form of
Dirac's conjecture on line arrangements of 1951\footnote{This conjecture
is different from the one addressed in \cite{tao}.} \cite{dirac}, which
will be stated in the next section.

In the following section we introduce several notions and facts used
throughout the paper. Section~\ref{approx} presents approximation
algorithms for GSL and analyzes their performance ratios. In particular,
we identify classes of line arrangements for which an approximate
solution with bounded performance ratio can be obtained. We also examine
structural relations between $g(M)$ and $w(M)$.

\section{Preliminaries}
\label{prelim}

\subsection{Some notions and facts}

The algorithms presented in this paper rely on well‑known procedures for solving basic problems of computational geometry: 
finding all intersection points of a set of lines, and determining all lines passing through each intersection point. 
Both problems can be solved by the $O((n+m)\log n)$ algorithm of Bentley and Ottmann \cite{bentleyOttmann}, 
the $O(m + n \log^2 n / \log \log n)$ algorithm of Chazelle \cite{chazelle}, 
and the optimal $O(m + n \log n)$ algorithm of Balaban \cite{balaban}. 
In what follows, we assume that Balaban’s algorithm is used. 
It can also be adapted to determine the intersection points lying on each line, see \cite{TCS}. 

Throughout, $M$ denotes a set of $n \ge 2$ lines in the plane, not all parallel.  
Given a graph or hypergraph $G$ (resp. a set of lines $M$), let $v$ be a vertex of $G$ (resp. an intersection point of lines from $M$).  
We denote by $d(v)$ the \emph{degree} of $v$, i.e., the number of edges incident to $v$ (resp. the number of lines of $M$ passing through $v$).  
The quantity $\Delta = \Delta(M)$ denotes the \emph{degree} of $M$, i.e., the maximum number of lines intersecting at a single point.

A \emph{pencil} of lines is a set of concurrent lines, and a \emph{near‑pencil} is a set of $n \ge 3$ lines such that $n-1$ of them are concurrent.  
By $conv(S)$ we denote the convex hull of a set $S$.  
An \emph{ordinary point} is an intersection point of exactly two lines.

\begin{fact}
A set of $n$ lines can contain at most $\frac{n(n-1)}{2}$ ordinary points, and this is the maximum possible number of intersections in a line arrangement. 
This bound is attained when every pair of lines intersects and no three lines pass through the same point.  
We call such a configuration an \emph{ordinary line arrangement}.
\label{fact1}
\end{fact}

A line in $M$ is called \emph{maximal} if it contains the maximum number $w(M)$ of intersections, and \emph{minimal} if it contains the minimum number of intersections, denoted $w'(M)$.

\begin{fact}
The values $w(M)$ and $w'(M)$ may coincide, as in an ordinary line arrangement where $w(M)=w'(M)=n-1$.  
They may also differ significantly; for example, in a near‑pencil, $w(M)=n-1$ while $w'(M)=2$.
\label{fact2}
\end{fact}

Let ${\cal A}$ be an approximation algorithm for a minimization problem $\Pi$.  
For an input $I \in \Pi$, let ${\cal A}(I)$ denote the solution returned by ${\cal A}$, and let $opt(I)$ denote the value of an optimal solution.  
The \emph{performance ratio} of ${\cal A}$ is defined as 
\[
R_{\cal A}(I) = \frac{|{\cal A}(I)|}{opt(I)}.
\]

We note that the minimum number of guards $g(M)$ can be arbitrarily smaller than $w(M)$.  
For example, if $M$ is a near‑pencil of $n$ lines, then $g(M)=2$ while $w(M)=n-1$.  
On the other hand, it is possible that $g(M)=w(M)$, e.g., when $M$ consists of several parallel lines intersected by a single transversal line.
\begin{fact}
If $M$ contains a line that is not parallel to any other line of $M$, then $g(M) \le w(M)$.
\label{fact3}
\end{fact}

This follows because the intersection points on such a line cover all lines in $M$.

\subsection{Guarding a set of lines}

Let $M$ be partitioned into classes $R_1, R_2, \dots, R_h$ such that two lines are parallel iff they belong to the same class.  
Let $r_i = |R_i|$ for $i=1,\dots,h$, and assume $r_1 \ge r_2 \ge \dots \ge r_h$.

Megiddo and Tamir \cite{megiddo} (with reference to a scheduling problem \cite[p.~54]{scheduling}) established the following fact.
\begin{fact}
If all intersection points of $M$ are ordinary, then Algorithm~1 finds an optimal guarding set of 
$\max\{ r_1, \lceil n/2 \rceil \}$ points in polynomial time.
\label{th1}
\end{fact}
\begin{algorithm}
\caption{: for the case where all intersection points are ordinary}
\textbf{Input:} Set $M$ of lines with no intersections of more than two lines.\\
\textbf{Output:} Minimum guarding set $S$ for the lines of $M$. 
\begin{algorithmic}
\State $S := \emptyset$;
\State \textbf{until} no lines remain in $M$ \textbf{or} no intersecting lines remain in $M$ \textbf{do}
\State Choose an arbitrary line $l_1$ from $R_1$ and an arbitrary line $l_2$ from $R_2$;
\State If $p = l_1 \cap l_2$ then set $S := S \cup \{p\}$;
\State Remove $l_1$ and $l_2$ from $M$;
\State Update the classes $R_i$ to maintain $r_i \ge r_j$ for $i < j$;
\State Return $M$ and $S$;
\State If $M = \emptyset$ then stop;
\State Otherwise, place a guard on each of the parallel lines remaining in $M$;
\State Let $S'$ be the set of these guards; then $S := S \cup S'$ is a guarding set for the original $M$.
\end{algorithmic}
\label{ordinary}
\end{algorithm}

Megiddo and Tamir also proved:

\begin{theorem}
GSL is strongly NP-hard.
\label{th2}
\end{theorem}

The problem is also known to be APX-hard, meaning that no polynomial-time algorithm can find a solution arbitrarily close to the optimum (unless P = NP).

\subsection{Dirac's conjecture}

Next we recall Dirac's conjecture \cite{dirac} and related results.

\begin{description}
\item[Dirac's Conjecture.] 
\emph{Given a finite set $A$ of non‑collinear points and the set $P(A)$ of all lines containing at least two points of $A$, 
there exists a point incident to at least $\lfloor n/2 \rfloor$ lines of $P(A)$.}
\end{description}

By line‑point duality, this conjecture can be reformulated as follows:\footnote{According to Gr\"unbaum \cite{grunbaum1}, most early discussions of the problem appeared in this dual form.}

\begin{description}
\item[Dirac's Conjecture, dual form.] 
\emph{Given a set of lines in the plane, not all parallel or concurrent, there exists a line containing at least $\lfloor n/2 \rfloor$ intersection points.}
\end{description}

Dirac proved a lower bound $\lceil \sqrt{n+1} \rceil$ \cite{dirac}, and for decades this remained the best known bound.  
Erd\H os proposed the following weaker form:

\begin{description}
\item[Weak Dirac's Conjecture.] 
\emph{Given a set of points (resp. lines) in the plane, there exists a point incident to at least $\lfloor n/c \rfloor$ lines (resp. a line containing at least $\lfloor n/c \rfloor$ intersections) for some positive constant $c$.}
\end{description}

This conjecture was proved by Beck \cite{beck} and Szemer\'edi–Trotter \cite{szemeredi_trotter} with unspecified $c$, later by Payne and Wood \cite{payne_wood} with $c=37$, and refined by Pham and Phi \cite{pham_phi} to $\frac{n}{26}+2$.

The original Dirac conjecture was shown to be false by several authors \cite{akiyama,felsner,grunbaum1,grunbaum_catalogue,grunbaum2,grunbaum3}.  
One related conjecture was recently proved by Han \cite{han}, yielding the following dual form:

\begin{theorem}
Every set $M$ of $n$ lines in the plane, which are non‑parallel and do not intersect at the same point, 
contains a line that has at least $\lceil n/3 \rceil + 1$ intersection points with other lines of $M$.
\label{th3}
\end{theorem}

This bound is attained for certain line arrangements and may therefore be best possible.

Note that it is possible to have $w(M)=\lceil n/3 \rceil + 1$ and this to be the minimum number of guards covering all lines of $M$, i.e., $w(M)=g(M)$.  
This holds, for example, when $M$ consists of three lines forming a triangle, or six lines consisting of a triangle and its medians.  
Other examples include the configurations $A(9,1)$ and $A(15,1)$ from Gr\"unbaum’s catalogue \cite{grunbaum_catalogue}.  
However, $g(M)$ may be significantly smaller than $w(M)$, e.g., for a near‑pencil of $n$ lines for large $n$.

\section{Approximation algorithms for GSL}
\label{approx}

\subsection{Set Cover Approximation}
\label{set cover}

The simplest heuristic solution for GSL is a set cover greedy strategy.
The GSL problem can be formulated as a set cover instance in which
the universe consists of the lines of $M$ and each intersection point
corresponds to the subset of lines passing through that point.
At each step an intersection of highest degree is selected and added to the guarding set.  
The algorithm terminates when all lines are covered.  
Consequently, the classical greedy set cover algorithm yields an
approximation ratio equal to the harmonic number
\[
H_n = 1+\frac12+\cdots+\frac1n,
\]
which satisfies
\[
H_n = \ln n + O(1).
\]
Hence GSL admits a polynomial-time
\(O(\log n)\)-approximation.

We now present a variant of that algorithm which better supports the performance ratio analysis. 
It is a combination of a set cover greedy procedure and Algorithm~\ref{ordinary}. 

\begin{algorithm} 
\caption{: set cover greedy + Megiddo--Tamir (hybrid approach)}
\textbf{Input:} Set $M$ of lines.\\
\textbf{Output:} Guarding set $S$ for the lines of $M$. 
\begin{algorithmic}
\State $i=0$;
\State $C_1 := \emptyset$;
\State {\em Stage 1:} \textbf{while} there exists an intersection of degree $\geq 3$ \textbf{do}
\State \quad Choose an intersection point $p$ with highest degree;
\State \quad $C_1 := C_1 \cup \{p\}$;
\State \quad Remove from $M$ all lines through $p$;
\State \quad Update degrees of remaining intersection points;
\State \quad $i := i+1$;
\State Return $C_1$ and the set $M'$ of remaining lines;
\State {\em Stage 2:} Apply Algorithm~\ref{ordinary} on $M'$ and obtain a guarding set $C_2$;
\State Compose $C = C_1 \cup C_2$ as a guarding set for $M$.
\end{algorithmic}
\label{mixed}
\end{algorithm}

Denote Algorithm~\ref{mixed} by $A_{\ref{mixed}}$.  
The approximation ratio of Algorithm~$A_{\ref{mixed}}$ is $\Omega(\log n)$. Indeed, the first stage of $A_{\ref{mixed}}$ coincides with the standard greedy set-cover algorithm applied to the line-cover instance. 
As shown by Hassin and Megiddo \cite{hassinMegiddo}, the greedy approach may exhibit a logarithmic worst-case performance ratio for line-covering problems. Since $A_{\ref{mixed}}$ 
makes exactly the same choices during Stage~1, it inherits this logarithmic lower bound. The second stage is invoked only after all intersections of degree at least three have been processed.
Thus it can only affect the lower-order part of the solution, does not change the asymptotic behavior, and therefore does not eliminate the asymptotic gap already accumulated during the greedy phase. 
Consequently, the worst-case approximation ratio of $A_{\ref{mixed}}$ remains $\Omega(\log n)$.


\subsubsection{Elementary bounds for the hybrid algorithm}

The following simple bounds on the performance ratio of Algorithm~\ref{mixed} will be used throughout Section~\ref{set cover}.  
Their proof is a routine case analysis of Stage~1 and is included for completeness.

\begin{proposition}
Let $A_{\ref{mixed}}$ be the hybrid algorithm. Then:
\[
\begin{aligned}
R_{A_{\ref{mixed}}}(M) &\le \frac{n}{2r_1} + \frac12, &
R_{A_{\ref{mixed}}}(M) &\le \frac{n}{3r_1} + 1, \\
R_{A_{\ref{mixed}}}(M) &\le \frac{\Delta}{2} + \frac12, &
R_{A_{\ref{mixed}}}(M) &\le \frac{\Delta}{3} + \frac{\Delta r_1}{n}.
\end{aligned}
\]
\label{th312}
\end{proposition}

\noindent {\bf Proof \ }
At each iteration of Stage~1, the algorithm selects an intersection point incident to at least three lines and removes all lines passing through that point.  
Thus Stage~1 removes at least three lines per iteration.

Three situations may occur:

\emph{(a) All lines are removed in Stage~1.}  
Then at most $n/3$ points are selected, giving  
\[
R \le \frac{n}{3r_1}, \qquad
R \le \frac{\Delta}{3}.
\]
\emph{(b) Only degree‑2 intersections remain.}  
The remaining lines form disjoint pairs, optimally guarded by at most $\lceil (n-3t)/2\rceil$ points, where $t$ is the number of Stage~1 iterations.  
Since $t \le n/3$, this yields  
\[
R \le \frac{n}{2r_1} + \frac12, \qquad
R \le \frac{\Delta}{2} + \frac12.
\]
\emph{(c) Only parallel lines remain.}  
Each remaining line contributes one guard.  
Since Stage~1 removes at least three lines per iteration, the number of remaining lines is at most $n - 3t$, giving  
\[
R \le \frac{n}{3r_1} + 1, \qquad
R \le \frac{\Delta}{3} + \frac{\Delta r_1}{n}.
\]
These bounds establish the stated inequalities.
\hfill$\square$

\subsubsection{Structural tools}
The following is a plain fact.

\begin{fact}
Let $U$ be a set of $m \geq 2$ points in the plane and let $d \geq 1$ lines pass through each of these points. 
Denote by $\tau_U(m,d)$ the minimum possible total number of lines passing through the given points over all possible line arrangements.
Let $\tau(m,d) = \min_U \tau_U(m,d)$. Then:
\begin{enumerate}
\item If $m_1 \leq m_2$, then $\tau(m_1,d) \leq \tau(m_2,d)$;
\item If $d_1 \leq d_2$, then $\tau(m,d_1) \leq \tau(m,d_2)$.
\end{enumerate}
\label{fact111}
\end{fact}

The following is a generalization of the well-known handshaking lemma to uniform hypergraphs \cite{berge}.

\begin{lemma}
Let $H(V,E)$ be an $r$-uniform hypergraph with $r \geq 2$.  
Then $e(H)= \frac{1}{r} \sum_{v \in V} d(v)$.
\label{lemma_hypergraph} 
\end{lemma}

\begin{corollary}
Let $H(V,E)$ be a hypergraph such that each hyperedge contains at most $r \geq 2$ vertices. 
Then $e(H) \geq \frac{1}{r} \sum_{v \in V} d(v)$.
\label{cor113}
\end{corollary}

\begin{definition}
Given a set ${\cal A}$ of points in the plane and a set ${\cal L}$ of lines through these points,  
define a hypergraph $H({\cal A},{\cal E})$ whose edges are
\[
{\cal E} = \{\, {\cal A} \cap l : l \in {\cal L} \,\}.
\]
\label{def}
\end{definition}

Let a maximal line of ${\cal A}$ contain $r$ points.  
Then Corollary~\ref{cor113} implies
\begin{equation}
|{\cal L}| \geq \frac{1}{r} \sum_{v \in {\cal A}} d(v).
\label{formula1}
\end{equation}

\begin{lemma}
In any arrangement of $n$ lines, the number of intersection points incident to $\Omega(n)$ lines is bounded by a constant.
\label{lemma99}
\end{lemma}

\noindent {\bf Proof \ }
Let $A$ be the set of intersection points of degree $\Omega(n)$.  
Assume for contradiction that $|A|=\omega(1).$\footnote{Recall that $f(n)=\omega(g(n))$ means that $f(n)$ ``dominates" $g(n)$ assimptotically. 
Formally, $f(n)=\omega(g(n))$ iff $\forall k>0 \ \exists n_0$, such that $\forall n>n_0$, $f(n) > kg(n)$.}
Let $L_A$ be the set of lines passing through points of $A$.

If no line contains two points of $A$, then $|L_A| = \omega(1)\cdot\Omega(n)=\omega(n)$, contradicting $|L_A|\le n$.

Thus some line contains at least two points of $A$.  
Consider the hypergraph $H(A,E)$ from Definition~\ref{def}.  
Let $r(n)$ be the maximum number of points of $A$ on any line.  
If $r(n)=\omega(1)$, then that line contributes $\omega(1)\cdot\Omega(n)=\omega(n)$ lines, again a contradiction.  
Hence $r(n)=O(1)$.

Let $d^*(n)=\Omega(n)$ be the minimum degree in $A$.  
Then by (\ref{formula1}),
\[
|E| \ge \frac{1}{r(n)} \sum_{v\in A} d(v) 
    \ge |A|\cdot d^*(n)
    = \omega(1)\cdot\Omega(n)
    = \omega(n).
\]
Since $|M|\ge |L_A|\ge |E|$, we obtain $|M|=\omega(n)$, a contradiction.  
Thus $|A|=O(1)$.
\hfill$\square$

\subsubsection{Constant-ratio classes}

\begin{proposition}
Let $P^*(M)$ be the set of intersection points incident with at least 3 lines.  
If $|P^*(M)| \leq c$ for some constant $c$, then Algorithm~\ref{mixed} solves GSL with performance ratio bounded by a constant. 
\label{th314}
\end{proposition}
\noindent {\bf Proof \ }
Let $p$ be the size of an optimal solution.  
Consider separately the three cases from Proposition~\ref{th312}.

(a) Then $|C_1|\le |P^*(M)|\le c$, so  
\[
R_{A_{\ref{mixed}}}\le \frac{c}{p}\le c.
\]
(If the instance is nontrivial (i.e. not all lines are concurrent), then $p \geq 2$, so  $R_{A_{\ref{mixed}}}\le \frac{c}{p}\le \frac{c}{2}$.)

(b) Let $q$ be the optimal solution size for $M'$.  
Then the algorithm finds at most $q+c$ points, and since $q\le p$,
\[
R_{A_{\ref{mixed}}}\le \frac{q+c}{p}\le 1+\frac{c}{p}\le 1+c.
\]

(c) Let $r$ be the number of parallel lines remaining.  
Then the algorithm finds at most $r+c$ points, giving
$$
R_{A_{\ref{mixed}}}\le \frac{r+c}{r_1}\le 1+\frac{c}{r_1}\le 1+c.
$$
\hfill$\square$

\medskip

Let $I$ be the set of intersection points of $M$, and suppose GSL has a solution $S$ with $|S|\le c$.  
Then one can test all subsets of $I$ of size at most $c$.  
Since $|I|\le {n\choose 2}=O(n^2)$, the number of tests is  
\[
{{n\choose 2}\choose c}=O(n^{2c}),
\]
a polynomial of degree $2c$.  
This is impractical for large $c$, and requires knowing $c$ in advance.  
Algorithm~\ref{mixed} avoids both issues.

\begin{proposition}
If the size of an optimal solution $S$ satisfies $|S|=\Omega(n)$ or $|S|=O(1)$, then Algorithm~\ref{mixed} finds a solution within a constant factor of the optimal. 
\label{prop_const_solution}
\end{proposition}
\noindent {\bf Proof \ }
First assume that $|S|=\Omega(n)$, that is,
\[
|S|=cn
\]
for some absolute constant $c$, $0<c<1$.
 Without loss of generality, suppose that
\[
r_1<\frac{n}{2}.
\]
Furthermore, the number $\alpha$ of steps performed during Stage~1 satisfies
\[
\alpha\leq \frac{n}{3}.
\]
Therefore, by Proposition~\ref{th312},
\[
R_{A_{\ref{mixed}}}
\leq
\frac{\alpha+r_1}{|S|}
\leq
\frac{\alpha+r_1}{cn}
<
\frac{\frac{n}{3}+\frac{n}{2}}{cn}
=
\frac{5}{6c}.
\]
Hence, in this case the approximation ratio is bounded by a constant.
 
\medskip
 
Now assume that $|S|=O(1)$, that is,
\[
|S|\leq c
\]
for some positive integer constant $c$.
 
Partition $S$ into three subsets:
\begin{itemize}
\item $S_1$, consisting of points incident to $\Omega(n)$ lines;
\item $S_2$, consisting of points incident to $\omega(1)$ and $o(n)$ lines;
\item $S_3$, consisting of points incident to $O(1)$ lines.
\end{itemize}
 
Since $|S|=O(1)$, we have
\[
|S_1|=O(1),\qquad |S_2|=O(1),\qquad |S_3|=O(1).
\]
 
Consider a point $p$ of degree $\Omega(n)$. If $p\notin S_1$, then all lines passing through $p$ must be covered by points from $S_2\cup S_3$. However, every point of $S_2\cup S_3$ is incident to only $o(n)$ lines, and since
$|S_2\cup S_3|=O(1)$, these points together can cover only $o(n)$ lines, which is insufficient to cover all $\Omega(n)$ lines through $p$. Hence every intersection point of degree $\Omega(n)$ must belong to $S_1$.
By the assumption of the theorem and Lemma~\ref{lemma99}, the number of such points is bounded by a constant.

Next consider a point $q$ of degree $\omega(1)$ and $o(n)$. After all points of degree $\Omega(n)$ have been selected, removing them decreases the degree of $q$ by at most a constant, since $|S_1|=O(1)$. Therefore $q$ remains of degree $\omega(1)$ and $o(n)$. If $q\notin S_2$, then its incident lines would have to be covered exclusively by points of $S_3$, each of which has degree $O(1)$. Since $|S_3|=O(1)$, this is impossible. Consequently, every intersection point of degree $\omega(1)$ and $o(n)$ must belong to $S_2$.
 
It follows that, before vertices of constant degree are considered, Algorithm~\ref{mixed} selects precisely the points of $S_1\cup S_2$. After their removal, all remaining vertices have degree $O(1)$. Since $|S_3|=O(1)$, the number of such remaining vertices is also bounded by a constant. 

It remains to consider the case in which all remaining uncovered lines become parallel during the construction of $S_3$.
In this case, an optimal completion is obtained by selecting one point on each remaining line. Since the number of such lines is bounded by a constant, this adds only a constant number of points and therefore does not affect the approximation analysis.
Thus, Stage~2 contributes only a constant number of points to the solution produced by Algorithm~\ref{mixed}.

Hence Algorithm~\ref{mixed} returns a guarding set of cardinality $O(1)$. Since the optimum solution satisfies $|S|=O(1)$, the approximation ratio is bounded by a constant as well.
  \hfill$\square$

\begin{theorem}
Let $M$ be a set of $n$ lines.  
Suppose that asymptotically the intersection points are either of degree $\Omega(n)$ or of bounded degree.  
Then Algorithm~\ref{mixed} has performance ratio bounded by a constant. 
\label{th315}
\end{theorem}

\noindent {\bf Proof \ }
Let $p$ be the size of an optimal solution $S$ of GSL on $M$.  
Let $A$ be the set of intersection points of degree $\Omega(n)$, and let $L_A$ be the set of all lines through points of $A$.

Consider first the case where an optimal solution does \emph{not} include some points from $A$.  
Then all lines through those points must be covered by other guards, so necessarily $p=\Omega(n)$.  
By Proposition~\ref{prop_const_solution}, when $|S|=\Omega(n)$ Algorithm~\ref{mixed} has performance ratio bounded by a constant.  

\sloppypar Otherwise, assume that $A$ is contained in every optimal solution.  
By Lemma~\ref{lemma99}, $|A|$ is bounded by a constant $\bar c$.  
In Stage~1, Algorithm~\ref{mixed} always selects an intersection of maximum degree, so it must eventually select all points of $A$ before selecting any point of bounded degree: 
after removing a constant number of lines through previously chosen points of $A$, the remaining points of $A$ still have degree $\Omega(n)$, strictly larger than the bounded degrees.

If, after selecting all points of $A$, no lines remain uncovered, then $A$ itself is a guarding set of bounded size.  
Hence $p$ is bounded, and the performance ratio of Algorithm~\ref{mixed} is bounded by a constant.

Otherwise, let $N$ be the set of lines that remain uncovered after all points of $A$ have been selected.  
Let $q$ be the size of an optimal solution $B^*$ for the minimum hitting set problem on $N$.  
Then any optimal solution $S$ for GSL on $M$ can be written as
\[
S = A \cup B^*, \qquad p = |S| = |A| + q.
\]
Since $N \subseteq M$, we have $q \le p$.  

By assumption, all intersection points of lines in $N$ are of bounded degree.  
Continuing Algorithm~\ref{mixed} on $N$ (Stage~1 on bounded‑degree intersections, followed by Algorithm~\ref{ordinary} on the remaining lines) yields an approximate solution $B$ of size $r$ such that
\[
\frac{r}{q} \le c'
\]
for some positive constant $c'$.

The algorithm’s final solution on $M$ is $A \cup B$, of size $|A| + r$.  
Thus the performance ratio satisfies
\[
R_{A_{\ref{mixed}}}
= \frac{|A| + r}{|A| + q}
\le \frac{\bar c + r}{q}
= \frac{\bar c}{q} + \frac{r}{q}
\le \bar c + c',
\]
which is a constant independent of $n$.  

In all cases, Algorithm~\ref{mixed} has performance ratio bounded by a constant, as claimed.
\hfill$\square$

\begin{remark}[Necessity of the conditions in Theorem~\ref{th315}]
If the intersection degrees are not restricted to the two regimes of Theorem~\ref{th315}, then the performance ratio of Algorithm~\ref{mixed} cannot be bounded by a constant.  
Indeed, the general set cover problem has worst‑case approximation ratio $\Omega(\log n)$, and GSL contains set cover instances as special cases (e.g., arrangements where each intersection corresponds to a set and each line corresponds to an element).  
Thus, unless P=NP, no polynomial‑time algorithm (including Algorithm~\ref{mixed}) can guarantee a constant approximation ratio for all line arrangements.  
This shows that the structural assumptions in Theorem~\ref{th315} are essential.
\end{remark}


\subsection{Relaxation technique}
We start with a statement that refines and rounds off those of Fact~\ref{th1} and Theorem~\ref{th2}.
It holds in the general case in which parallel lines as well as intersections of three or more lines may exist, and is a base for an approximation algorithm.  
\begin{theorem}
\begin{enumerate}
\item
If $r_1 \geq \left\lceil \frac{n}{2} \right\rceil$, then
$g(M)=r_1\leq w(M)$, and an optimal solution can be found in
polynomial time.
\item
If $r_1 < \left\lceil \frac{n}{2} \right\rceil$, then
\[
r_1 \le g(M)\le \left\lceil \frac{n}{2} \right\rceil
<\frac{3}{2}w(M),
\]
and the computation of an optimal solution is strongly NP-hard.
\end{enumerate}
\label{prop1}
\end{theorem}
\noindent {\bf Proof \ }
1.
Assume that
\[
r_1 \geq \left\lceil \frac{n}{2} \right\rceil ,
\]
and let
\[
R_1=\{a_1,a_2,\ldots,a_{r_1}\}.
\]
Let $b_1,b_2,\dots,b_s$, where $s=n-r_1$, be the lines not belonging to $R_1$.
Every such line intersects every line of $R_1$.
Since the lines of $R_1$ are pairwise disjoint, each of them must
contain a guard. Hence
\[
g(M)\ge r_1.
\]
Place guards at $a_1 \cap b_1, a_2 \cap b_2, \dots, a_s \cap b_s$, then place a guard at each of the lines $a_{s+1}, \dots, a_{r_1}$. 
Thus, all lines are guarded by $r_1$ guards, which is clearly  optimal. 
The last inequality holds since each line which is not in $R_1$ contains at least $r_1$ intersections with the lines of $R_1$,
and possibly other intersections with lines not in $R_1$.
An optimal solution can be found in polynomial time by the Algorithm~\ref{ordinary}.   

\smallskip

2. Now let $r_1 < \lceil \frac{n}{2} \rceil$. 
The first inequality holds since each line from $R_1$ must contain a guard. 
Towards proving the second inequality, we construct a new set of lines $M'$ with $|M'|=|M|$, 
which has classes of parallel lines $R_1',R_2',\dots,R_h'$, as $|R_i'|=|R_i|$ for any $i$, $1 \leq i \leq h$, and such that no three or more lines intersect at the same point. 
For this we perform the following procedure. 
\begin{description}
\item{(a)}
Choose a line $l$ which contains at least one intersection point of $l$ with at least two other lines. 
Let $d$ be the minimum distance of an intersection point not on $l$ to $l$. 
Replace $l$ by a line $l'$ which is parallel to $l$ and is at a distance $d/2$ from $l$.  
\item{(b)}
After finitely many repetitions of Step (a), we obtain a set of lines
$M'$ having the same parallel classes as $M$ and containing no
intersection point incident with three or more lines. Every multiple
intersection of $M$ is thus replaced by a collection of ordinary
intersections in $M'$.
\end{description}

Let $x$ be a point of intersection of at least $t \geq 3$ lines of $M$.
In $M'$ these lines are replaced by a set of other $t$ lines. 
Denote it by $F_x$.
Fact~\ref{fact1} implies that in $M'$ instead of the single intersection point $x$ there will be ${t \choose 2}$ intersection points of pairs of lines from $F_x$.
An optimal solution $S$ to the GSL problem on $M’$ is found in polynomial time by Algorithm \ref{ordinary}  and counts $\lceil \frac{n}{2} \rceil$ intersection points.



Let ${\cal F}(S)$ denote the family of all sets $F_x$ for which
$S$ contains at least one of the corresponding ordinary intersection
points.

Define
\[
S_1=\{x:\;F_x\in{\cal F}(S)\},
\]
and let
\[
S_2
=
S-
\{\text{all points of }S\text{ belonging to some }F_x\in{\cal F}(S)\}.
\]

Then $S_1\cup S_2$ is a guarding set for $M$ with at most $\lceil \frac{n}{2} \rceil$ elements and thus $g(M) \leq \lceil \frac{n}{2} \rceil$.
Indeed, each point $x\in S_1$ guards all lines represented by the family
$F_x$, while all remaining guarded lines are still covered by points of
$S_2$.

Moreover,
\[
|S_1\cup S_2|\le |S|,
\]
since several points of $S$ associated with the same family $F_x$ are
replaced by the single point $x$.


\smallskip

The third inequality is implied by Theorem~\ref{th3}: since $w(M) \geq  \lceil \frac{n}{3} \rceil + 1$, we have
$$
\frac{3}{2}w(M) \geq \frac{3}{2}\left(\left\lceil \frac{n}{3} \right\rceil + 1\right) \geq \frac{3}{2}\left(\frac{n}{3}+1\right) = \frac{n}{2} + \frac{3}{2} > \left\lceil \frac{n}{2} \right\rceil. 
$$ 
The strong NP-hardness of the problem when $r_1 < \lfloor \frac{n}{2} \rfloor$ follows from \cite{megiddo}.
\ \ \ \ \ \ \ \ \ \ \ \ \hfill$\square$

\bigskip

Regarding maximal lines, we have:
\begin{itemize}
\item
If $r_1 > \lceil \frac{n}{2} \rceil$, then only the lines not from $R_1$ can be maximal and $w(M) \geq r_1 > \lceil \frac{n}{2} \rceil$.
\item
If $r_1 = \lceil \frac{n}{2} \rceil$, then there are maximal lines not from $R_1$; some (or all) lines from $R_1$ can be maximal, too. 
\end{itemize}
By the first part of Theorem~\ref{prop1}, in either case $g(M)=r_1 \leq w(M)$, and an optimal solution can be found in polynomial time as demonstrated in the proof.
In the nontrivial case $r_1 < \lceil \frac{n}{2} \rceil$, by the second part of Theorem~\ref{prop1} we have  $g(M) \leq \lceil \frac{n}{2} \rceil$. 
In fact, this is the case making the GSL problem NP-hard.

Theorem \ref{prop1} shows that the structural dichotomy known for ordinary arrangements extends to arbitrary line arrangements containing parallel lines and multiple intersections. The perturbation argument provides a mechanism for transferring solutions between the degenerate and nondegenerate settings.
The theorem  provides a clear structural boundary between the tractable and intractable regimes of the guarding-lines problem. When the dominant parallel class is sufficiently large, the geometry forces every remaining line to intersect it, yielding an optimal solution of size $r_1$ computable in polynomial time. In contrast, when the dominant class is smaller, the problem becomes strongly NP-hard, yet the perturbation argument used in the proof shows that every instance admits a relaxed form with only ordinary intersections and an optimal guarding set of size at most $\lceil n/2 \rceil$. 
Such a guarding set can be found by the following algorithm with polynomial time complexity. 

\medskip

\begin{algorithm} 
\caption{: relaxation technique}
\textbf{Input:} Set $M$ of lines.\\
\textbf{Output:} Guarding set $S$ for the lines of $M$.
\begin{algorithmic}
\State Construct a set $M'$ of ordinary lines as in the proof of Theorem~\ref{prop1};
\State Using Algorithm~\ref{ordinary}, find an optimal solution $S$ of GSL on $M'$;
\State Identify the family ${\cal F}(S)$ of sets of lines $F_x$ which contain at least one element of $S$; 
\State Determine an approximate solution $S_{{\cal A}_1} = S_1 \cup S_2$, where:
\begin{itemize}
\item[]
$S_1 = \{ x : F_x \in {\cal F}(S) \}$;
\item[]
$S_2 = S - \{ x : F_x \in {\cal F}(S)  \}$.
\end{itemize}
\end{algorithmic}
\label{relaxation}
\end{algorithm}

In view of Fact~\ref{fact1} and elementary analytic geometry concerning line intersections, it is easy to realize that 
the algorithm finds an approximation solution with a polynomial number of arithmetic operations, and that the numbers produced in the course of 
the computation are bounded by a polynomial in the line coefficients.   

Denote for short Algorithm \ref{relaxation} by $A_{\ref{relaxation}}$. 
Note that if some sets $F_y$ contribute to $A_{\ref{relaxation}}(M)$ with $k>1$ elements, then $|A_{\ref{relaxation}}(M)|$ will be strictly smaller than $\lceil \frac{n}{2} \rceil$, 
as the value of  $|A_{\ref{relaxation}}(M)|$ can be obtained by reducing $\lceil \frac{n}{2} \rceil$ by $k-1$ for each such $F_y$. 
Otherwise, we would have $|A_{\ref{relaxation}}(M)| = \lceil \frac{n}{2} \rceil$. 
About the performance ratio  of Algorithm \ref{relaxation} we have the following straightforward estimation.



\begin{proposition}
$R_{A_{\ref{relaxation}}}(M) \leq \frac{1}{2} \min \{ \Delta, \frac{n}{r_1} \}+1$.
\label{alg}
\end{proposition}
\noindent {\bf Proof \ }
First notice that $opt(M)  \geq \frac{n}{\Delta}$, $opt(M) \geq r_1$, and recall that $\Delta \leq n-1$, $r_1 \geq 1$, and $|A_{\ref{relaxation}}(M)| \leq \lceil \frac{n}{2} \rceil$. 
Then we have
\begin{align*}
&\frac{|A_{\ref{relaxation}}(M)|}{opt(M)} \leq \frac{\lceil \frac{n}{2} \rceil}{n/\Delta} < \frac{(\frac{n}{2}+1)\Delta}{n} = \\
&=\frac{\Delta}{2}+\frac{\Delta}{n}\leq \frac{\Delta}{2}+\frac{n-1}{n}<\frac{\Delta}{2}+1,
\end{align*}
and 
\begin{align*}
&\frac{|A_{\ref{relaxation}}(M)|}{opt(M)} \leq \frac{\lceil \frac{n}{2} \rceil}{r_1} < \frac{(\frac{n}{2}+1)}{r_1} =\\
&= \frac{n}{2r_1}+\frac{1}{r_1} \leq \frac{n}{2r_1}+1.
\end{align*}
\hfill$\square$

\begin{corollary}
\begin{description}
\item
If $\Delta \leq c_1$ for some constant $c_1$, then $R_{A_{\ref{relaxation}}}(M) \leq \frac{c_1}{2}+1$.
\item
If $r_1 \geq \frac{n}{c_2}$ for some constant $c_2$, then $R_{A_{\ref{relaxation}}}(M) \leq \frac{c_2}{2}+1$.
\end{description}
\label{cor_alg}
\end{corollary}

Although Proposition~\ref{alg} does not yield a logarithmic worst-case
approximation ratio, the relaxation procedure exploits geometric
features of the arrangement that are ignored by purely combinatorial
covering algorithms. It therefore appears to be a promising framework
for obtaining stronger approximation guarantees on restricted classes
of line arrangements, a question left for future investigation.
One can also first use the relaxation directly to obtain a simple $\lceil n/2 \rceil$-bounded approximation. 
Then, by interpreting the relaxed instance as a set-cover problem over a geometrically meaningful candidate family, one can derive a hybrid algorithm that combines geometric preprocessing with the classical greedy set-cover strategy.


\subsection{Approximation using $l_{\max}$ and $w(M)$}

Let $M$ be a set of $n$ lines in the Euclidean plane. For a point
$x$ lying on a line $l\in M$, let $B_x$ denote the set of all lines of $M$
passing through $x$. We call $x$ the \emph{root} of the {\em beam} $B_x$.

The following general framework underlies several approximation algorithms
for the GSL problem.

\paragraph{General Scheme}
\begin{enumerate}
\item Choose a line $l\in M$ together with its intersection points;
\item For each line $l'$ parallel to  $l$, place a guard at the intersection point of $l'$ and another line intersecting $l'$;
\item For every beam with a root on $l$ whose lines are not yet completely guarded, place a guard at its root.
\end{enumerate}

\subsubsection*{Algorithm Based on a Minimal Line}

Let $w'(M)$ denote the minimum number of intersection points lying on a line
of $M$.

\begin{algorithm}
\caption{Algorithm $A_{\min}$}
\textbf{Input:} Set $M$ of lines.\\
\textbf{Output:} Guarding set $S$.
\begin{algorithmic}[1]
\State Identify a line $l_{\min}$ containing $w'(M)$ intersection points;
\State Place guards at all intersection points on $l_{\min}$;
\State For every line parallel to $l_{\min}$, place a guard at an arbitrary
      intersection point on that line.
\end{algorithmic}
\end{algorithm}

Algorithm $A_{\min}$  is clearly polynomial. 
Its performance ratio satisfies
\[
R_{A_{\min}}  = \frac{w'(M)+r_1}{g(M)} \leq \frac{w'(M)+r_1}{r_1} = \frac{w'(M)}{r_1}+1.
\]
Thus $A_{\min}$ achieves a bounded ratio whenever $w'(M)/r_1=O(1)$.
It also has a `small' performance ratio for any line arrangement $M$ with a `small' $w'(M)$.   
For example, for a near pencil with $n$ lines 
the performance ratio of the Algorithm~$A_{\min}$ equals $\frac{w'(M)}{r_1} + 1 = \frac{2}{1}+1=3$. 


\subsubsection*{Greedy Algorithm Based on a Maximal Line}

Let $l_{\max}$ be a line containing $w(M)$ intersection points. Let
\[
x_1,x_2,\ldots,x_t
\]
be the intersection points on $l_{\max}$ ordered so that
\[
|B_{x_1}|\leq |B_{x_2}|\leq \cdots \leq |B_{x_t}|.
\]
\begin{algorithm}
\caption{Algorithm $A_{\max}$}
\textbf{Input:} Set $M$ of lines.\\
\textbf{Output:} Guarding set $S$.
\begin{algorithmic}[1]
\State Compute $l_{\max}$ and the beams rooted at its intersection points.
\State Process the beams in nondecreasing order of size.
\State Whenever possible, use lines from the current beam to guard uncovered
      lines parallel to $l_{\max}$.
\State After all parallel lines have been covered, place a guard at the root
      of every beam containing an unguarded line.
\end{algorithmic}
\end{algorithm}
\begin{fact}
\label{fact:greedy-optimal}
Within the class of algorithms that following the General Scheme, Algorithm
$A_{\max}$ uses the minimum possible number of guards.
\end{fact}

\noindent {\bf Proof \ }
Let $p$ be the number of intersection points on $l_{\max}$ and let $r$ be the
number of lines parallel to $l_{\max}$. Stage~2 of the General Scheme uses
lines from beams to guard parallel lines. Every beam whose lines are used in
Stage~2 reduces the number of root guards needed in Stage~3 by exactly one.

Let $z$ be the number of beams whose lines are fully used in Stage~2. Then
the total number of guards placed is
\[
r + (p - z).
\]
To minimize this quantity, one must maximize $z$. A beam contributes to $z$
only if \emph{all} of its lines are used in Stage~2. Using a beam of size~$s$
requires guarding $s$ parallel lines. Thus smaller beams are easier to use
fully. Processing beams in nondecreasing order of size maximizes the number
of beams that can be fully used before parallel lines are exhausted. Hence
$A_{\max}$ minimizes the number of guards among all algorithms following the
General Scheme.
\hfill$\square$

\begin{fact}
\label{fact:beam-small}
Suppose that at most $2i$ lines intersect a line $l$ at $i$ points. Then at
least $\lfloor i/2\rfloor$ of the corresponding beams have cardinality at
most $2$.
\end{fact}
\noindent {\bf Proof \ }
Let the beams rooted at the $i$ intersection points have sizes
$s_1,\ldots,s_i$. Then
\[
s_1+\cdots+s_i\le 2i.
\]
Assume fewer than $\lfloor i/2\rfloor$ beams have size at most~$2$. Then at
least $\lceil i/2\rceil+1$ beams have size at least~$3$. The total number of
lines intersecting $l$ is then at least
\[
3\Bigl(\Bigl\lceil\frac{i}{2}\Bigr\rceil+1\Bigr)
+1\cdot\Bigl(\Bigl\lfloor\frac{i}{2}\Bigr\rfloor-1\Bigr)
=i+2\Bigl\lceil\frac{i}{2}\Bigr\rceil+2.
\]
Since $\lceil i/2\rceil\ge i/2$, this quantity exceeds $2i$, contradicting
the assumption. Thus at least $\lfloor i/2\rfloor$ beams have size at most~2.
\hfill$\square$

By Theorem~\ref{th3}, if $n=3k$ then $w(M) \geq k+1$, and if $n=3k+1$ or $n=3k+2$ then $w(M) \geq k+2$. 
In what follows we consider more in detail the first case, the others being analogous. 
We have the following theorem.
\begin{theorem}
\label{thm:main-structure}
Let $w(M)=k+m$, where $1\le m\le 2k-1$, and suppose that $l_{\max}$ belongs
to a parallel class consisting of $j$ lines.

\begin{enumerate}
\item If $m\ge \lfloor k/2\rfloor$, then $g(M)\le w(M)$.
\item If $m<\lfloor k/2\rfloor$, then
\[
g(M)\le w(M)+\Bigl\lfloor\frac{j-1}{2}\Bigr\rfloor.
\]
\end{enumerate}
In both cases the guard locations can be computed in polynomial time.
\end{theorem}
\noindent {\bf Proof \ }
We analyze the number of lines intersecting $l_{\max}$ and the distribution
of beam sizes.

\medskip\noindent\textbf{Case 1: $m\ge\lfloor k/2\rfloor$.}
There are $3k-j$ lines intersecting $l_{\max}$ at $k+m$ points. Consider the
quantity
\[
((k+m)-(j-1))\cdot 2 + (j-1)\cdot 1.
\]
The first term counts the contribution of intersection points incident with
two nonparallel lines; the second counts those incident with one. The
inequality
\[
((k+m)-(j-1))\cdot 2 + (j-1)\cdot 1 \ge 3k-j
\]
is equivalent to $m\ge k/2 - 1/2$, which is implied by $m\ge\lfloor k/2\rfloor$.

Thus at least $j-1$ intersection points on $l_{\max}$ are incident with
exactly one nonparallel line. These $j-1$ lines can guard the $j-1$ parallel
lines distinct from $l_{\max}$. Placing guards at the remaining
$(k+m)-(j-1)$ intersection points yields a guarding set of size $k+m=w(M)$.

\medskip

\noindent\textbf{Case 2: $m < \lfloor k/2 \rfloor$.}

Let $l_{\max}$ be one of the $j$ parallel lines, where $1 \leq j \leq k+m$. Then, once again, $3k-j$ lines intersect $l_{\max}$ at $k+m$ points.

We have
\begin{equation}
((k+m)-(j-1))\cdot 3 + (j-2)\cdot 2 + 1\cdot 1 \geq 3k-j,
\label{ineq_case2}
\end{equation}
since the above inequality is equivalent to $m\geq 0$, which is always true.

Inequality~(\ref{ineq_case2}) implies that $l_{\max}$ contains at least one intersection point $p$ that is incident to exactly one intersecting line, and that there are at least $j-1$ intersection points (including $p$) that are incident to at most two intersecting lines.

The extremal case requiring the largest number of guards occurs when $p$ is the only point whose intersection multiplicity is one. In this case, there are at least $j-2$ points whose intersection multiplicity is two.

Let $g$ be the line intersecting $l_{\max}$ at $p$. 
Place a guard at the intersection of $g$ and $l$, thereby guarding both lines.
Now let $q$ be an intersection point on $l_{\max}$ of multiplicity two, and let $g_1$ and $g_2$ be the two lines intersecting $l_{\max}$ at $q$. 
Place a guard at the intersection of $g_1$ with an uncovered line parallel to $l_{\max}$, 
and another guard at the intersection of $g_2$ with a different uncovered line parallel to $l_{\max}$. In this way, four lines are guarded.

Repeat this procedure until all $j-1$ lines parallel to $l_{\max}$ are covered. 
In either case, exactly one guard is placed on each line, for a total of $j-1$ guards. 
The guards are located at intersections with lines belonging to $\Bigl\lfloor \frac{j-1}{2} \Bigr\rfloor + 1$ beams rooted on $l_{\max}$: one consisting of a single line intersecting $l_{\max}$, while each of the remaining beams consist of two such lines).

Next we place one guard at each of the remaining
\[
(k+m)-\Bigl\lfloor \frac{j-1}{2} \Bigr\rfloor-1
\]
intersection points on $l_{\max}$.
Thus, the total number of guards required to cover all lines is
\begin{align*}
&\left((k+m)-\Bigl\lfloor \frac{j-1}{2} \Bigr\rfloor-1\right)\cdot 1
+\Bigl\lfloor \frac{j-1}{2} \Bigr\rfloor\cdot 2+1\cdot 1 \\
&=(k+m)+\Bigl\lfloor \frac{j-1}{2} \Bigr\rfloor =w(M)+\Bigl\lfloor \frac{j-1}{2} \Bigr\rfloor.
\end{align*}
Clearly, if $l_{\max}$ contains more than one intersection point of multiplicity one, then the total number of guards required to cover all lines can be reduced. 

\smallskip

When $j=1$ (i.e., when no other line is parallel to $l_{\max}$), the bound reduces to the obvious inequality $g(M) \leq w(M)$, stated in Fact~\ref{fact3}. If $j=k+m$, the upper bound becomes
\begin{align*}
&(k+m)+\Bigl\lfloor \frac{k+m-1}{2} \Bigr\rfloor \\
&\leq (k+m)+\frac{k+m-1}{2} \\
&< \frac{3}{2}(k+m)  = \frac{3}{2}w(M).
\end{align*}
\hfill$\square$

The guard placement configurations for both cases can be found efficiently. 
Computing the intersection points along $l_{\max}$ and tracking their line multiplicities requires a single pass over the line system, taking $O(n^2)$ time. 
The subsequent greedy pairing of parallel lines in Case~2 and the allocation of baseline guards on $l_{\max}$ take at most $O(n^2)$ time. 
Because the algorithm only needs to store the line coordinates and an intersection matrix, the framework requires $O(n^2)$ memory. 
Thus, the guard configuration is determined in strictly $O(n^2)$ time and space.

\medskip

For `small' or `large' values of $j$, certain refinements or alternative estimates for the number of required guards are possible.

\begin{proposition}
\label{prop:small-j}
If $j\le \lfloor k/2\rfloor$, then
\[
g(M)<\frac{5}{4}w(M).
\]
\end{proposition}
\noindent {\bf Proof \ }
Since $j\le\lfloor k/2\rfloor$, we have
\[
\Bigl\lfloor\frac{j}{2}\Bigr\rfloor
\le \Bigl\lfloor\frac{\lfloor k/2\rfloor}{2}\Bigr\rfloor
\le \Bigl\lfloor\frac{k}{4}\Bigr\rfloor
\le \frac{k}{4}.
\]
Since $w(M)=k+m\ge k+1$, we have $w(M)>k$ and thus
\[
\frac{k}{4}<\frac{w(M)}{4}.
\]
By Theorem~\ref{thm:main-structure},
\begin{align*}
&g(M)\le w(M)+\Bigl\lfloor\frac{j-1}{2}\Bigr\rfloor <  w(M)+\Bigl\lfloor\frac{j}{2}\Bigr\rfloor \leq \\
&w(M)+\Bigl\lfloor\frac{k}{4}\Bigr\rfloor < w(M)+\frac{w(M)}{4} =\frac{5}{4}w(M).
\end{align*}
\hfill$\square$


\begin{proposition}
\label{prop:large-j}
Let $j=k+i$, where $1\le i\le m<\lfloor k/2\rfloor$. Then
\[
g(M)\le \frac{3}{2}w(M)-i+1.
\]
\end{proposition}
\noindent {\bf Proof \ }
There are $2k-i$ lines intersecting $l_{\max}$ at $k+m$ points. Consider the
quantity
\[
(k-i)\cdot 2 + (m+i)\cdot 1.
\]
This counts the total number of lines intersecting $l_{\max}$ under the
assumption that $m+i$ beams have size~1 and the remaining $k-i$ beams have
size~2. The inequality
\[
(k-i)\cdot 2 + (m+i)\cdot 1 \ge 2k-i
\]
holds for all $m\ge 0$, so at least $m+i$ beams have size~1.

These beams guard $m+i$ parallel lines. The remaining number of intersecting
lines is
\[
(2k-i)-(m+i)=2k-m-2i,
\]
and the remaining number of beams is
\[
(k+m)-(m+i)=k-i.
\]
Since $2k-m-2i<2(k-i)$, Fact~\ref{fact:beam-small} implies that at least
$\lfloor (k-i)/2\rfloor$ of these beams have size at most~2. These beams
guard the remaining $k-i$ parallel lines.

The number of unused beams is
\[
(k-i)-\Bigl\lfloor\frac{k-i}{2}\Bigr\rfloor
\le \frac{k+i}{2}-i+1.
\]
Placing guards at their roots yields a total of
\begin{align*}
&\left(\frac{k+i}{2}-i+1\right)+2\Bigl\lfloor\frac{k-i}{2}\Bigr\rfloor+(m+i) \\
&\le \frac{3}{2}(k+m)-i+1=\frac{3}{2}w(M)-i+1.
\end{align*}
\hfill$\square$

\begin{remark}
Note that the condition $i \leq m \leq \left \lfloor \frac{k-1}{2} \right\rfloor$ used in the proof is not a restriction of the generality. 
In fact, $i > m$ would imply that $j=k+i>k+m=w(M)$. 
Then every line intersecting a maximal line would contain strictly more intersection points than that maximal line, contradicting its maximality.
Moreover, if $m \geq  \lfloor \frac{k}{2} \rfloor$, then Theorem~\ref{thm:main-structure} implies that $w(M)$ guards suffices to guard all lines of set $M$.  
\label{remark12}
\end{remark}
\begin{proposition}
\label{prop:ratio-maxline}
\begin{enumerate}
\item If $w(M)<n/2$, then
\[
R_{A_{\max}}<\frac{1}{2}\min\!\left\{
\Delta\Bigl(\frac{r_1}{n}+1\Bigr),\,
\frac{n}{r_1}+1
\right\}.
\]
\item If $w(M)\ge n/2$, then
\[
R_{A_{\max}}\le \min\!\left\{
\frac{w(M)}{r_1},\,
\frac{w(M)\Delta}{n}
\right\}.
\]
\end{enumerate}
\end{proposition}
\noindent {\bf Proof \ }
By Theorem~\ref{thm:main-structure},
\[
R_{A_{\max}}\le \frac{w(M)+\lfloor j/2\rfloor}{r_1}
\le \frac{w(M)+r_1/2}{r_1}
=\frac{w(M)}{r_1}+\frac{1}{2}.
\]
Similarly,
\[
R_{A_{\max}}\le \frac{(w(M)+r_1/2)\Delta}{n}.
\]
If $w(M)<n/2$, then
\[
R_{A_{\max}} \le  \frac{w(M)}{r_1}+\frac{1}{2}
<\frac{n}{2r_1}+\frac{1}{2}
=\frac{1}{2}\left(\frac{n}{r_1}+1\right),
\]
and
\[
R_{A_{\max}} \le  \frac{(w(M)+r_1/2)\Delta}{n}
<\frac{(n/2+r_1/2)\Delta}{n}
=\frac{\Delta}{2}\left(\frac{r_1}{n}+1\right).
\]
If $w(M)\ge n/2$, then Theorem~\ref{thm:main-structure} implies that
$A_{\max}$ uses at most $w(M)$ guards, giving the stated bounds.
\hfill$\square$

\medskip

We define a line arrangement $M$ as \emph{non-trivial} if it excludes configurations that are exactly solvable in polynomial time, such as arrangements where $g(M) > n/2$ due to large sets of mutually parallel lines, or pencil arrangements.

\begin{theorem}
\label{thm:final}
Let $M$ be a non-trivial arrangement of $n$ lines. Then:
\begin{enumerate}
\item If $g(M) = \Omega(w(M))$, Algorithm $A_{\max}$ achieves a constant approximation ratio.
\item If $g(M) = o(w(M))$, then $M$ contains both subclasses that admit a constant approximation ratio and subclasses whose hardness matches the general APX-hardness of the problem.
\end{enumerate}
\end{theorem}
\noindent {\bf Proof \ }
1. Suppose $g(M) \ge c \cdot w(M)$ for some constant $c > 0$. The approximation ratio $R_{A_{\max}}$ satisfies
\begin{align*}
&R_{A_{\max}} \leq  \frac{w(M)+r_1}{g(M)} \leq   \frac{w(M)+r_1}{c \cdot w(M)} \leq \\              
&\frac{(n-1) + \lceil \frac{n}{2} \rceil}{c (\lceil \frac{n}{3} \rceil +1)} <  
\frac{(n-1) + \frac{n}{2} +1}{c (\frac{n}{3} +1)} = \\
&\frac{3n/2}{c (n/3+1)}  = 
\frac{9n}{2cn+6c} < \frac{9}{2c},
\end{align*}
confirming that Algorithm $A_{\max}$ guarantees a constant-factor approximation.
For non-trivial arrangements, $g(M) \le n/2$ and $w(M) \ge n/3 + 1$, which bounds the constant $c \le 3/2$. 
Thus, the analytical upper bound $\frac{9}{2c}$ evaluates to at least $3$.

\medskip
  
2. Consider an arrangement with $n$ a perfect square, in which $\sqrt{n}$ points are selected such that exactly $\sqrt{n}$ lines pass through each selected point, and no other multiple intersections occur.
For this class of instances, $g(M) = \sqrt{n}$ while $w(M) \ge n/3 + 1$, yielding $g(M) = o(w(M))$. 
On these instances, the greedy algorithm routinely finds the optimal solution, yielding a perfect approximation ratio of 1. 

Conversely, consider subclasses derived from standard geometric reductions of NP-hard hitting set instances. When configured such that $g(M) = o(w(M))$, 
these arrangements preserve the underlying hardness gap, meaning their approximability remains strictly bounded away from 1 by the problem's inherent APX-hardness.
\hfill$\square$



\section{Lines and Segments}

All questions pertinent to sets of straight lines can also be asked about sets of straight line segments. 
Various aspects of the problem of {\em guarding a set of segments} (GSS) were considered in \cite{compMath,TCS,DAM}.    
In particular, it was shown that GSS is strongly NP-hard \cite{TCS} and approximation algorithms for the problem were proposed and scrutinized theoretically and experimentally. 
Next we exhibit a plain relation between GSS and GSL. 

Let $\Pi_{L}$ and $\Pi_{S}$ be an arbitrary search problem that can be defined on a set $M$ of lines and line segments, respectively.
Suppose that a solution to an individual problem $I$ may be a set of intersection points of lines/segments, and/or
a set of segments whose end-points are such intersections.
The intersections and the segments formed by $M$ may possibly have labelings, which may be involved in the solution. 
In case of a decision problem the solution may be just {\em Yes} or {\em No}. 
In the above setting, we have the following fact. 

\begin{proposition}
If $\Pi_{L}$ is NP-hard, then $\Pi_{S}$ is NP-hard. 
If $\Pi_{S}$ is in $P$, then $\Pi_{L}$ is in $P$. 
\label{th33}
\end{proposition}
\noindent {\bf Proof \ }
Let $\Pi_{L}$ be an NP-hard problem defined on a set $M$ of lines in the plane. 
Let $P(M)$ be the set of all intersection points of lines from $M$ and $S(M)$ the set of all segments with end-points from $P(M)$. 
Let $\hat M = M \cap conv (P(M))$, i.e., $\hat M$ is the set of segments obtained by taking the portion of each line which falls within $conv (P(M))$.
Then each segment from $S(M)$ is identical or is a part of some segment in $\hat M$. 
Hence, any subset of $P(M)$ and of $S(M)$ is contained in the line segments in $\hat M$.
Thus, any solution to $\Pi_{S}$ which may be a subset of $P(M)$ and/or $S(M)$ would be a solution to $\Pi_{L}$, and vice versa.
This establishes a polynomial reduction of $\Pi_{S}$ to $\Pi_{L}$, which completes the proof of both parts of the proposition. 
\hfill$\square$

\medskip

Proposition~\ref{th33}, coupled with Theorem~\ref{th2}, implies:

\begin{corollary}
GSS is strongly NP-hard.
\label{cor33}
\end{corollary}

Proposition~\ref{th33} says that NP-hardness results on a set of lines may imply NP-hardness results on a set of line segments, 
since, because of the exhibited polynomial reduction, a problem defined on a set of lines can be regarded as a particular case of a problem defined on a set of line segments. 
By that reason, structural results and polynomial algorithms that hold for problems defined on a set of line segments hold for problems defined on a set of lines, as well.   
On occasions, the opposite may hold, too. 
For example, all algorithms presented in this paper can easily be adapted for GSS.  


\section{Concluding remarks}
\label{conclusion}
In this paper we investigated the problem of guarding a set of lines
in the plane and studied several approximation approaches based on
set-cover techniques, geometric relaxations, and structural properties
of line arrangements. The obtained results reveal that the performance
of approximation algorithms is closely related to incidence parameters
such as the maximum intersection multiplicity and the number of
intersections on a maximal line.

Several questions remain open. In particular, it would be desirable
to determine whether logarithmic approximation guarantees can be
improved for broader classes of arrangements, and whether stronger
relations between the guarding number and classical incidence
parameters can be established. Another direction is the extension of
the methods developed here to higher-dimensional arrangements and to
other geometric covering problems.

We conclude by stating the following  conjecture.
\begin{conjecture}
$g(M) \leq w(M)$.
\end{conjecture}

\section*{Acknowledgments}
The author thanks J\'anos Pach for pointing out a useful reference and for helpful comments at an early stage of this research.

\end{document}